\documentclass[sn-mathphys-num]{sn-jnl}

\usepackage{graphicx}%
\usepackage{multirow}%
\usepackage{amsmath,amssymb,amsfonts}%
\usepackage{amsthm}%
\usepackage{mathrsfs}%
\usepackage[title]{appendix}%
\usepackage{xcolor}%
\usepackage{textcomp}%
\usepackage{manyfoot}%
\usepackage{booktabs}%
\usepackage{algorithm}%
\usepackage{algorithmicx}%
\usepackage{algpseudocode}%
\usepackage{listings}%

\theoremstyle{thmstyleone}%
\newtheorem{theorem}{Theorem}
\newtheorem{corollary}[theorem]{Corollary}
\newtheorem{lemma}[theorem]{Lemma}%

\theoremstyle{thmstyletwo}%
\newtheorem{example}{Example}%
\newtheorem{remark}{Remark}%

\theoremstyle{thmstylethree}%
\newtheorem{definition}{Definition}%

\begin{document}

\title[An Order Relation between Eigenvalues and Symplectic Eigenvalues]{An Order Relation between Eigenvalues and Symplectic Eigenvalues of a Class of Infinite-Dimensional Operators}


\author[1]{\fnm{Tiju Cherian} \sur{John}}\email{tijucherian@fulbrightmail.org}
	\equalcont{These authors contributed equally to this work.}
	
	\author[2]{\fnm{V. B. Kiran} \sur{Kumar}}\email{vbk@cusat.ac.in}
	\equalcont{These authors contributed equally to this work.}
	
	\author*[3]{\fnm{Anmary} \sur{Tonny}}\email{anmarytonny97@gmail.com}

	\affil[1]{\orgdiv{Wyant College of Optical Sciences, 1630 E University Blvd, The University of Arizona, Tucson, Arizona 85721}}
	
	\affil[2]{\orgdiv{Department of Mathematics, Cochin University of Science and Technology, Kerala, India 682022}}
	
	\affil*[3]{\orgdiv{Department of Mathematics, Cochin University of Science and Technology, Kerala, India 682022}}	


\abstract{In this article, we obtain some results in the direction of ``infinite-dimensional symplectic spectral theory". We prove an inequality between the eigenvalues and the symplectic eigenvalues of a special class of infinite-dimensional operators. 
 Let $T$ be a positive invertible operator such that $T - \alpha I$ is compact for some $\alpha > 0$.  
 Denote by $\{{\lambda_j^R}^\downarrow(T)\}$,  the set of eigenvalues of $T$ lying strictly to the right side of $\alpha$, arranged in decreasing order and let $\{{\lambda_j^L}^\uparrow(T)\}$ denote the set of eigenvalues  of $T$ lying strictly to the left side of $\alpha$, arranged in increasing order. Furthermore, let $\{{d_j^R}^\downarrow(T)\}$ denote the symplectic eigenvalues of $T$ lying strictly to the right of $\alpha$, arranged in decreasing order and $\{{d_j^L}^\uparrow(T)\}$ denote the set of symplectic eigenvalues of $T$ lying strictly to the left of $\alpha$, arranged in increasing order, respectively (such an arrangement is possible as it is shown in the article that the only possible accumulation point for the symplectic eigenvalues is $\alpha$). Then we show that 
		$${d_j^R}^\downarrow(T) \leq {\lambda_j^R}^\downarrow(T), \quad j = 1,2, \dots , s_R$$ and  $${\lambda_j^L}^\uparrow(T) \leq {d_j^L}^\uparrow(T), \quad j = 1,2, \dots , s_L,$$ where $s_R$ and $s_L$ denote the number of symplectic eigenvalues of $T$ lying strictly to the right and left of $\alpha$, respectively. 
		This generalizes a finite-dimensional result obtained by Bhatia and Jain (J. Math. Phys. 56, 112201 (2015)). The class of Gaussian Covariance Operators (GCO) and positive Absolutely Norm attaining Operators ($(\mathcal{AN})_+$ operators) appears as special cases of the set of operators that we consider.}

\keywords{Symplectic transformations; Symplectic spectrum; Interlacing Theorems; Gaussian Covariance Operators; Positive $\mathcal{AN}$ Operators.}

\pacs[MSC Classification]{Primary 47A10, 47B07, 81S10  ; Secondary 46L60, 47B15, 47B65}



\maketitle

\section{Introduction}
	Williamson's Normal form for $2n\times 2n$ real positive matrices is a symplectic analogue of the spectral theorem for normal matrices. With the recent developments in quantum information theory,  Williamson's normal form has opened up an active research area that may be dubbed as ``finite-dimensional symplectic spectral theory", analogous to the usual spectral theory and matrix analysis.  
 \par A symplectic matrix is a real matrix $M$ of order $2n,$ satisfying the identity $M^{T}JM = J$, where
	$ J = \begin{bmatrix}
		0 & I_n \\
		-I_n & 0
	\end{bmatrix}$ and $I_n
	$ is the identity matrix of order $n$.
	Williamson's normal form is a diagonalization result for positive definite matrices of even order, where the diagonalizing matrix comes from the symplectic group \cite{jw01}. A formal statement is given below:
	\begin{theorem} (Williamson's normal form \cite{jw01}) \label{finitewil}
		\newline Let $A$ be a positive definite real matrix of order $2n$. Then there exists a symplectic matrix $M$ and an $n \times n$ strictly positive diagonal matrix $D$ such that  
		\begin{equation}\label{wili}
			A = M^T
			\begin{bmatrix}
				D & 0 \\
				0 & D
			\end{bmatrix}M.
		\end{equation}	
		Furthermore, the diagonal matrix $D$ is unique up to the ordering of its entries.
	\end{theorem}
 
	\begin{definition}The diagonal entries of the matrix $D$ appearing in \eqref{wili} are called the \emph{symplectic eigenvalues} of $A$.
	\end{definition}	
	Symplectic eigenvalues play an important role in continuous variable quantum information theory \cite{AdRaSa14, mag01}. It is an invariant for quantum Gaussian states, and several properties of these states depend on the symplectic eigenvalues of their covariance matrices \cite{weedbrooketal2012, Par13}. Several researchers, including  Bhatia and Jain \cite{rjt01, bhatia_jain_2021};  Idel, Gaona and Wolf \cite{wolf}; Hiai and Lim  \cite{logmaj};  Jain and Mishra \cite{jain2022derivatives}; and  Son and Stykel \cite{tracemin} have made significant developments in studying the properties of symplectic eigenvalues in the finite-dimensional setting. These developments include the introduction of log-majorisations for symplectic eigenvalues, perturbation bounds for the normal form, derivatives and Lidskii type theorems for symplectic eigenvalues, and trace minimization theorems.	
	
	An infinite-dimensional analogue of Williamson's Normal form is obtained in  \cite{bv01} and has been very useful in the study of infinite-mode quantum Gaussian states. The $2n\times 2n$ real strictly positive matrix in the finite-dimensional Williamson's normal form (Theorem \ref{finitewil}) is replaced with a positive, invertible operator on $\mathcal{H} \oplus \mathcal{H}$, where $\mathcal{H}$ is a real separable Hilbert space. More details are available in \cite{bv01}.	
	Let $I$ be the identity operator on $\mathcal{H}$. The involution operator $J$ on $\mathcal{H} \oplus \mathcal{H}$ is defined by $ J = 
	\begin{bmatrix}
		0 & I \\
		-I & 0
	\end{bmatrix}$, where $I$ is the identity operator on $\mathcal{H}$.
	A bounded invertible linear operator $M$ on $\mathcal{H} \oplus \mathcal{H}$ is called a symplectic transformation if $M^{T}JM = J$. Now we state the infinite-dimensional analogue of the Williamson's normal form. 
	
	\begin{theorem} (Williamson's normal form \cite{bv01}) \label{rajatij}
		Let $\mathcal{H}$ be a real separable Hilbert space and $T$ be a positive invertible operator on $\mathcal{H} \oplus \mathcal{H}$. Then there exists a positive invertible operator $D$ on $\mathcal{H}$ and a symplectic transformation $M : \mathcal{H} \oplus \mathcal{H} \rightarrow \mathcal{H} \oplus \mathcal{H}$ such that 
		\[T = M^T
		\begin{bmatrix}
			D & 0 \\
			0 & D
		\end{bmatrix}M.\]
	\end{theorem} 
	\begin{remark} 		The operator $D$ in Theorem \ref{rajatij} is unique up to conjugation with an orthogonal transformation \cite{bv01}. The spectrum of $D$ is defined as the \textbf{symplectic spectrum} of $T$ and is denoted by $\sigma_{sy}(T)$. Similarly, the eigenvalues of the operator $D$ is defined as the symplectic eigenvalues of the operator $T$. 
 In this article, we are 
  concerned with operators that have a countable symplectic spectrum.
	\end{remark}

	Analysis of the symplectic spectrum in the infinite-dimensional setting is important in quantum information with infinite-mode states, and also in the theory of infinite-mode quantum Gaussian states \cite{john2018infinite}. But this subject is relatively new and there are only a  few results available in this direction. Nevertheless, a Szeg\H{o} type theorem for symplectic eigenvalues was proved  by Bhatia, Jain and Sengupta in \cite{bhatia2021szegHo}. In this article, we present an order relation between the eigenvalues and the symplectic eigenvalues of a class of infinite-dimensional operators under the setting of Theorem \ref{rajatij}. 

 This article is organized as follows. For the rest of this section, we discuss the main result of our article and also describe the preliminary results that we use in this work.In Section~\ref{mr},  after setting up the necessary background, we prove the main result of this article.  In Section~\ref{sec:special-cases}, we illustrate the main result on the classes of GCOs and $(\mathcal{AN})_+$ operators. The application of our main result to the case of GCO (Corollary~\ref{decinc}) is of independent interest.
	
	\subsection{Discussion on the Main Result} \label{discussionmainresults}
	The class of operators we consider in this article are 
	\textbf{strictly positive translations of compact operators} on $\mathcal{H}\oplus \mathcal{H}$ (where $\mathcal{H}$ is an infinite-dimensional real separable Hilbert space), that is, positive invertible operators $T$ such that $T-\alpha I$ is compact for some $\alpha>0$. This is a rich class of operators that contain the Gaussian Covariance Operators (GCO) and the positive absolutely norm attaining $((\mathcal{AN})_+)$
	operators. The inclusions between these operators are strict, and is explained in more detail in Appendix \ref{appendixinclusion}. 
	
	Now we turn to the statement of our main result. Let $\mathcal{H}$ be an infinite-dimensional real separable Hilbert space and let $T$ be a  positive invertible operator on  $\mathcal{H} \oplus \mathcal{H}$ as in Theorem \ref{rajatij}. Assume that there exists an $\alpha>0$ such that $T-\alpha I$ is a compact operator. Then $T$ has at most countably many eigenvalues with $\alpha$ as the only possible accumulation point. Hence, $T$ has either finitely many eigenvalues or $\alpha$ is an accumulation point of the eigenvalues of $T$. Let $\{{\lambda_j^R}^\downarrow(T)\}$ denote  the set of eigenvalues of $T$ lying strictly to the right side of $\alpha$ arranged in decreasing order and let $\{{\lambda_j^L}^\uparrow(T)\}$ denote the set of eigenvalues of $T$ lying strictly to the left side of $\alpha$ arranged in increasing order. We show in Lemma \ref{lemmacardinality} that the symplectic spectrum of $T$ is at most countable. Also, we show that $\alpha$ is the only possible accumulation point for the symplectic spectrum of $T$. Assuming this preliminary result, let $\{{d_j^R}^\downarrow(T)\}$ denote the set of symplectic eigenvalues of $T$ lying strictly to the right of $\alpha$, arranged in decreasing order and $\{{d_j^L}^\uparrow(T)\}$ denote the set of symplectic eigenvalues of $T$ lying strictly to the left of $\alpha$, arranged in increasing order. Let $s_R$ and $s_L$ denote the number of symplectic eigenvalues of $T$ lying strictly to the right and left of $\alpha$, respectively. The main result of this article (Theorem \ref{decincgen}) states that 
	 $${d_j^R}^\downarrow(T) \leq {\lambda_j^R}^\downarrow(T), \quad j = 1,2, \dots , s_R$$ and  $${\lambda_j^L}^\uparrow(T) \leq {d_j^L}^\uparrow(T), \quad j = 1,2, \dots , s_L.$$

	Our main result is Theorem \ref{decincgen}, which generalizes the finite-dimensional counterpart proved by Bhatia and Jain \cite[Theorem 11(ii)]{rjt01}. We illustrate this in Remark \ref{finitedimcounterpart}. The major tool used in the finite-dimensional context is the min-max principle. The min-max principle holds for the infinite-dimensional positive compact operators, but this technique cannot be applied to the operators considered in this article because they lack compactness. 
	Also, a careful ordering of the eigenvalues and symplectic eigenvalues is essential in our setting. We overcome these challenges by using the truncation method and an approximation result obtained in \cite{botcher}. The details are given in Section \ref{prelims}.
	
	\subsubsection{Preliminary Results} \label{prelims}
	
	In this section, we discuss some preliminaries needed to prove the main result. Let $\mathcal{K}$ be a complex separable Hilbert space, and $A \in B(\mathcal{K})$ be a self-adjoint operator. Define  $$m = \underset{\|x \| = 1}{\textrm{ inf }} \langle Ax, x  \rangle; \quad \text{and} \quad M = \underset{\|x \| = 1}{\textrm{ sup }} \langle Ax, x  \rangle.$$ 
	   Let $\nu$ and $\mu$ denote, respectively,  the minimum and maximum of the essential spectrum of $A$.  The eigenvalues lying outside $[\nu, \mu]$ are at most countable, that is, there are at most countable eigenvalues in the intervals $[m, \nu)$ and $(\mu, M]$ (see \cite{gohberg} for details). Therefore, the eigenvalues lying in the interval $[m, \nu)$ can be arranged in increasing order, and if they are countably infinite, they converge to $\nu$. Similarly, the eigenvalues in the interval $(\mu, M]$ can be arranged in  decreasing order, and if they are countably infinite, they converge to $\mu$. Now we arrange the eigenvalues in the interval $[m, \nu)$ as $${\lambda_1^L}^\uparrow(A) \leq {\lambda_2^L}^\uparrow(A) \leq \cdots \leq {\lambda_S^L}^\uparrow(A),$$ and the eigenvalues in the interval $(\mu, M]$ as $${\lambda_1^R}^\downarrow(A) \geq {\lambda_2^R}^\downarrow(A) \geq \cdots \geq {\lambda_Q^R}^\downarrow(A),$$(counting multiplicities), where $S,Q \in \{1,2,\dots \} \cup \{ \infty \}$. 
   
   Now let us fix an orthonormal basis $\{e_1, e_2, \dots \}$  of the separable Hilbert space $\mathcal{K}$. For $n = 1,2, \dots$, define $P_n$ to be the orthogonal projection of $\mathcal{K}$ onto the span of $\{ e_1, e_2, \dots , e_n\}$. Put $A_n = P_nAP_n$, which can be considered as the truncation of $A$ and can be identified with an $n \times n$ matrix. Denote the eigenvalues of $A_n$ (counting multiplicities) by $$\lambda_1(A_n), \lambda_2(A_n), \dots, \lambda_n(A_n).$$ Now for $j = 1,2, \dots, n$, denote the eigenvalues by $\lambda_j^\uparrow(A_n)$, when arranged in increasing order and the eigenvalues by $\lambda_j^\downarrow(A_n)$, when arranged in decreasing order. A well-known result in the Galerkin method of spectral approximation \cite[Theorem 3.1]{botcher} serves as an important tool in our proof techniques. A formal statement of this result adapted to a form, that is, useful for us is provided in Theorem \ref{mnnresult}. Essentially, this result states that, for each $k = 1,2, \dots$, as $n \rightarrow \infty$, the sequence of eigenvalues $\{\lambda_k^\uparrow(A_n)\}$ converges to ${\lambda_k^L}^\uparrow(A)$ (if there are countably infinite eigenvalues for $A$ to the left of $\nu$) or converges to $\nu$ (if there are only a finite number of eigenvalues for $A$ to the left of $\nu$). A similar statement holds for the eigenvalues lying to the right of $\mu$ as well.
	The theorem also concludes that as $k,n \rightarrow \infty$, the sequences $\{\lambda_k^\uparrow(A_n)\}$ and $\{\lambda_k^\downarrow(A_n)\}$ converge to $\nu$ and $\mu$, respectively. Now we give the formal statement of this result.
	
	\begin{theorem}\cite[Theorem 3.1]{botcher} \label{mnnresult}
		Let $A \in B(\mathcal{K})$ be self-adjoint. Then $$\underset{n \rightarrow \infty}{\lim} \lambda_k^\downarrow(A_n) = \begin{cases}
			{\lambda_k^R}^\downarrow(A), \textrm{ if } Q = \infty, \textrm{ or } 1 \leq k \leq Q \textrm{ when $Q<\infty$}, \\
			\mu, \textrm{ if } Q < \infty \textrm{ and } k \geq Q + 1,
		\end{cases}$$
		$$\underset{n \rightarrow \infty}{\lim} \lambda_k^\uparrow(A_n) = \begin{cases}
			{\lambda_k^L}^\uparrow(A),  \textrm{ if } S = \infty, \textrm{ or } 1 \leq k \leq S \textrm{ when $S < \infty$}, \\
			\nu, \textrm{ if } S < \infty \textrm{ and } k \geq S + 1,
		\end{cases}$$ In particular, 
		$$\underset{k \rightarrow \infty}{\lim} \underset{n \rightarrow \infty}{\lim} \lambda_k^\downarrow(A_n) = \mu, \quad \underset{k \rightarrow \infty}{\lim} \underset{n \rightarrow \infty}{\lim} \lambda_k^\uparrow(A_n) = \nu.$$
	\end{theorem} 

Now we  discuss some basic facts about bounded linear operators on a real Hilbert space $\mathcal{H}$. First, we define the complexification of a real Hilbert space and the complexification of a real linear operator defined on the real Hilbert space.

	\begin{definition} \label{defncomplexification} \cite{bv01}
        By the \emph{complexification} of a real Hilbert space $\mathcal{H}$ we mean the complex Hilbert space $\hat{\mathcal{H}} = \mathcal{H} + i \cdot \mathcal{H} = \{x + i \cdot y : x,y \in \mathcal{H} \}$ with addition, complex-scalar product and inner product defined in the obvious way. For a bounded operator $A$ on the real Hilbert space $\mathcal{H}$, define an operator $\hat{A}$ on the complexification $\hat{\mathcal{H}}$ of $\mathcal{H}$ by $\hat{A}(x + i \cdot y) = Ax + i \cdot Ay$. The operator $\hat{A}$ thus defined is called the \emph{complexification} of $A$. Define the spectrum of $A$, denoted by $\sigma(A)$, to be the spectrum of $\hat{A}$. Thus the eigenspectrum of $A$ and $\hat{A}$ are the same. The transpose of $A$ is the unique operator $A^T$ that satisfies the equation $\left\langle A^Tx,y \right\rangle_{\mathcal{H}}=\left\langle x,Ay \right\rangle_{\mathcal{H}}$, for all $x,y\in \mathcal{H}$. We say that the operator $A$ is self-adjoint if $A=A^T$.
    \end{definition}

    
    The lemma below states some basic properties of the complexified operator. These come in handy later in this article.

    \begin{lemma} \label{lemmacomplexification} 
    Let $\mathcal{H}$ be a real Hilbert space.
        \begin{enumerate}
			\item \label{item:1.1} If $A$ is a compact operator on $\mathcal{H}$, then $\hat{A}$ is compact on $\hat{\mathcal{H}}$. 
			
			\item \label{item:2.1} For operators $A_1$ and $A_2$ on $\mathcal{H}$, $\widehat{A_1A_2} = \hat{A_1}\hat{A_2}$. 
			
			\item \label{item:3.1} For any operator $A$ on $\mathcal{H}$, the adjoint of $\hat{A}$, $(\hat{A})^* = \widehat{(A^T)}.$ Thus the operator $\hat{A}$ is Hermitian on $\hat{\mathcal{H}}$ if and only if the operator $A$ is self-adjoint on $\mathcal{H}$.

            \item \label{item:4.1}  For any invertible operator $A$ on $\mathcal{H}$, $\hat{A}^{-1} = \widehat{A^{-1}}.$
   
			\item \label{item:5.1} For any positive operator $A$ on $\mathcal{H}$, $\widehat{\sqrt{A}} = \sqrt{\hat{A}}.$  
		\end{enumerate}
    \end{lemma}

    \begin{proof}
        The proofs of \ref{item:1.1}, \ref{item:2.1}, \ref{item:3.1}, and \ref{item:4.1} are straightforward. 
        \newline \underline{Proof of \ref{item:5.1}:} Since $A$ is a positive operator, it is self-adjoint and hence, $\hat{A}$ is Hermitian on $\hat{\mathcal{H}}$. Now since the spectrum of $A$ and $\hat{A}$ coincides, $\hat{A}$ is a positive operator on $\hat{\mathcal{H}}$. 
        Now the result follows from part \ref{item:2.1} of the Lemma.
    \end{proof}
    
    The following result is implicit in the proof of \cite[Theorem 4.3]{bv01}, but we state and prove it here for the sake of completeness.
	\begin{lemma} \label{lemmasymspecdefn}
	    Let $\mathcal{H}$ be a real separable Hilbert space and $A$ be a positive invertible operator on $\mathcal{H} \oplus \mathcal{H}$. Then the symplectic spectrum of $A$ coincides with the positive part of the spectrum of the operator $i\widehat{\sqrt{A}J\sqrt{A}}$.
	\end{lemma}

    \begin{proof}
        Consider the operator $B = \sqrt{A}J\sqrt{A}$ which is a skew-symmetric operator on $\mathcal{H} \oplus \mathcal{H}$. From \cite[Theorem 3.4]{bv01}, $B$ has a decomposition
        $$
        B = V^T \begin{bmatrix}
            0 & -P \\
            P & 0
        \end{bmatrix}V,
        $$ where $V$ is an orthogonal operator on $\mathcal{H} \oplus \mathcal{H}$ and $P$ is a positive invertible operator on $\mathcal{H}$. Hence, the spectrum of $B$,  $\sigma(B) = \sigma\left(V^T \begin{bmatrix}
            0 & -P \\
            P & 0
        \end{bmatrix}V\right) = \sigma\left( \begin{bmatrix}
            0 & -P \\
            P & 0
        \end{bmatrix} \right)$. Therefore, $\sigma(B) = i\sigma(P) \cup - i \sigma(P)$. That is,
        \begin{align*}
            -\sigma(P) \cup \sigma(P) &= i \sigma(B) = i \sigma(\hat{B}) 
            = \sigma(i \hat{B}) \\ &= \sigma( i \widehat{\sqrt{A}J\sqrt{A}}).
        \end{align*}
        Since $P$ is a positive operator, the  equation above gives that the spectrum $\sigma(P)$ coincides with the positive part of the spectrum of the operator $i \widehat{\sqrt{A}J\sqrt{A}}$.
        Also from \cite[Theorem 4.3]{bv01} $A$ has the Williamson's normal form
        $$
        A = L^T \begin{bmatrix}
            P & 0 \\
            0 & P
        \end{bmatrix}L$$ for the symplectic transformation $L$ on $\mathcal{H} \oplus \mathcal{H}$ and the same positive operator $P$ on $\mathcal{H}$. Hence, by definition, we have that the symplectic spectrum of $A$, $\sigma_{sy}(A) = \sigma(P)$. Therefore, it follows that the symplectic spectrum of $A$ coincides with the positive part of the spectrum of the operator $i\widehat{\sqrt{A}J\sqrt{A}}$.
    \end{proof}

    \begin{remark}
        It follows from the Lemma above that the symplectic eigenvalues of a positive invertible operator $A$ coincides with the positive part of the eigenvalues of the operator $i\widehat{\sqrt{A}J\sqrt{A}}$.
    \end{remark}
 
	\section{The Main Result} \label{mr}
	In this section, we state and prove the main result (Theorem \ref{decincgen}) of this article. 
	Throughout this section, \textbf{$\mathcal{H}$ denotes an infinite-dimensional real separable Hilbert space}. 
	Before proceeding to the main result, we make some observations and also prove a few lemmas about operators on a real Hilbert space and their complexifications. These results are used in the proof of the main theorem. Recall the operator $J$ defined before Theorem \ref{rajatij}. Note that this operator satisfies the properties $J = -J^T$ and $J^2 = -I$. First, we prove two lemmas about positive invertible operators $T$ such that $T -\alpha I$ is compact for some $\alpha>0$.


	\begin{lemma} \label{lemmaoncompactness}
		Let $\mathcal{H}$ be a real separable Hilbert space, and $T$ be a positive invertible operator on $\mathcal{H} \oplus \mathcal{H}$. If $T - \alpha I$ is compact for some $\alpha > 0$, then the operator $(\sqrt{T}J\sqrt{T})^T(\sqrt{T}J\sqrt{T}) - \alpha^2 I$ is compact. Moreover if $T - \alpha I$ is of finite rank for some $\alpha > 0$, then the operator $(\sqrt{T}J\sqrt{T})^T(\sqrt{T}J\sqrt{T}) - \alpha^2 I$ is also of finite rank.          
	\end{lemma}

	\begin{proof}
		We prove the lemma in two steps. First, we show that the compactness of the operator $T - \alpha I$ implies the compactness of the operator $(JT)^2 + \alpha^2 I$. Secondly, we show that the compactness of the operator $(JT)^2 + \alpha^2 I$ is equivalent to the compactness of the operator $(\sqrt{T}J\sqrt{T})^T(\sqrt{T}J\sqrt{T}) - \alpha^2 I$.
        Let $L = T - \alpha I$. Then $L$ is a compact operator by our assumption. Now
        \begin{align*}
   (JT)^2 + \alpha^2 I&= J(L + \alpha I)J(L + \alpha I) + \alpha^2 I\\ 
   &=(JL+\alpha J)(JL+\alpha J)+\alpha^2 I\\
   &=(JL)^2 + \alpha JLJ - \alpha L + \alpha^2J^2 + \alpha^2 I\\
   &=(JL)^2 + \alpha JLJ - \alpha L + \alpha^2(-I) + \alpha^2 I.\\
   & = (JL)^2 + \alpha JLJ - \alpha L.
\end{align*} The compactness of $L$ implies that the operator on the right side of the  equation above is compact. Hence, the operator $(JT)^2 + \alpha^2 I$ is compact. 
		Now to prove the conclusion of our statement, note that \begin{align*}
			(JT)^2 + \alpha^2 I \textrm{ is compact } &\Leftrightarrow JTJT + \alpha^2 I \textrm{ is compact} \\
			&\Leftrightarrow \sqrt{T}JTJT(\sqrt{T})^{-1} + \alpha^2 I \textrm{ is compact} \\
            &\Leftrightarrow \sqrt{T}JTJ\sqrt{T}\sqrt{T}(\sqrt{T})^{-1} + \alpha^2 I \textrm{ is compact} \\
			&\Leftrightarrow \sqrt{T}J\sqrt{T}\sqrt{T}J\sqrt{T} + \alpha^2 I \textrm{ is compact} \\
			&\Leftrightarrow -(\sqrt{T}J\sqrt{T})^T(\sqrt{T}J\sqrt{T}) + \alpha^2 I \textrm{ is compact} \\
			&\Leftrightarrow (\sqrt{T}J\sqrt{T})^T(\sqrt{T}J\sqrt{T}) - \alpha^2 I \textrm{ is compact.}
		\end{align*} Hence, the operator $(\sqrt{T}J\sqrt{T})^T(\sqrt{T}J\sqrt{T}) - \alpha^2 I$ is compact when $T - \alpha I$ is compact. The same proof holds with compactness replaced by the finite rank property.
	\end{proof}

	\begin{lemma} \label{lemmacardinality}
		Let $T$ be a positive invertible operator on $\mathcal{H} \oplus \mathcal{H}$, where $\mathcal{H}$ is a real separable Hilbert space. If $T - \alpha I$ is compact for some $\alpha > 0$, then the symplectic spectrum of $T$ is at most countable.
	\end{lemma}
	
	\begin{proof}
		 From Lemma \ref{lemmaoncompactness}, we have that the operator $(\sqrt{T}J\sqrt{T})^T(\sqrt{T}J\sqrt{T}) - \alpha^2 I$ is compact, hence, the operator $\widehat{(\sqrt{T}J\sqrt{T})^T(\sqrt{T}J\sqrt{T})} - \alpha^2 I$ is compact. That is, $\widehat{(\sqrt{T}J\sqrt{T})^T(\sqrt{T}J\sqrt{T})}$ has a countable spectrum. Now
		$$\widehat{(\sqrt{T}J\sqrt{T})^T(\sqrt{T}J\sqrt{T})} = \widehat{(\sqrt{T}J\sqrt{T})^T}\widehat{(\sqrt{T}J\sqrt{T})} = -i\widehat{(\sqrt{T}J\sqrt{T})^T}i\widehat{(\sqrt{T}J\sqrt{T})}.$$ Hence, the operator $i\widehat{(\sqrt{T}J\sqrt{T})^T}i\widehat{(\sqrt{T}J\sqrt{T})}$ has a countable spectrum. Note that the spectrum of $i\widehat{(\sqrt{T}J\sqrt{T})^T}i\widehat{(\sqrt{T}J\sqrt{T})}$ is the square of the spectrum of $\widehat{i\sqrt{T}J\sqrt{T}}$. Therefore, $i\widehat{\sqrt{T}J\sqrt{T}}$ has a countable spectrum. But from Lemma \ref{lemmasymspecdefn}, the symplectic spectrum of the operator $T$ coincides with the positive part of the spectrum of the operator $i\widehat{\sqrt{T}J\sqrt{T}}$.  Hence, the symplectic spectrum of $T$ is at most countable. 
	\end{proof}
	
	Our main result (Theorem \ref{decincgen}) provides a relationship between the eigenvalues and the symplectic eigenvalues of a positive invertible operator $T$ on the separable Hilbert space $\mathcal{H} \oplus \mathcal{H}$ such that $T - \alpha I$ is compact for some $\alpha > 0$. This makes $\alpha$ to be the only possible accumulation point for the spectrum of $T$. In this situation, either $T$ has finitely many eigenvalues or $\alpha$ is an accumulation point of the eigenvalues of $T$. Let $\{{\lambda_j^R}^\downarrow(T)\}$ denote  the set of eigenvalues of $T$ lying strictly to the right side of $\alpha$, arranged in the decreasing order and  $\{{\lambda_j^L}^\uparrow(T)\}$ denote the set of eigenvalues of $T$ lying strictly to the left side of $\alpha$, arranged in the increasing order. Also, the essential spectrum of the operator $i\widehat{\sqrt{T}J\sqrt{T}}$ is $\{-\alpha, \alpha\}$. 
 This can be seen as follows. Define $M = i\widehat{\sqrt{T}J\sqrt{T}}$. Then $M$ is self-adjoint and  
	\begin{align*}
		0 \leq M^2 &= (i\widehat{\sqrt{T}J\sqrt{T}})^2 
		= (i\widehat{\sqrt{T}J\sqrt{T}})(i\widehat{\sqrt{T}J\sqrt{T}}) \\
		&= -(\widehat{\sqrt{T}J\sqrt{T}})(\widehat{\sqrt{T}J\sqrt{T}}) 
		= (\widehat{\sqrt{T}J\sqrt{T}})^T(\widehat{\sqrt{T}J\sqrt{T}}) \\
		&= (\widehat{\sqrt{T}J\sqrt{T}})^T(\widehat{\sqrt{T}J\sqrt{T}}) - \alpha^2I + \alpha^2I \\
		M^2 &= K + \alpha^2I, 
	\end{align*} where $K = (\widehat{\sqrt{T}J\sqrt{T}})^T(\widehat{\sqrt{T}J\sqrt{T}}) - \alpha^2I$. Hence, $K$ is compact by Lemma \ref{lemmaoncompactness}. Hence, $\alpha ^2$ is the only element in the essential spectrum of $M^2$. Hence, the essential spectrum of $M$ is at most $\{- \alpha, \alpha\}$. Since the operator $i\widehat{\sqrt{T}J\sqrt{T}}$ is self-adjoint (hence normal), the essential spectrum of $M$ is $\{- \alpha, \alpha \}$ (Corollary 3.3 of \cite{bv01}). Therefore the eigenvalues of $i\widehat{\sqrt{T}J\sqrt{T}}$ accumulate around $\alpha$ and $-\alpha$, symmetrically about $0$ on the real axis (when they are countably infinite). Hence, the set of symplectic eigenvalues of $T$ accumulates about $\alpha$ when they are countably infinite. So the set of symplectic eigenvalues of $T$ is either finite or accumulates at $\alpha$. Let $\{{d_j^R}^\downarrow(T)\}$ denote  the set of symplectic eigenvalues of $T$ lying strictly to the right side of $\alpha$, arranged in decreasing order and let $\{{d_j^L}^\uparrow(T)\}$ denote the set of symplectic eigenvalues of $T$ lying strictly to the left side of $\alpha$, arranged in increasing order. 
	
	Now we develop some prerequisites required to prove our main result. Let $\{e_1, e_2, \dots \}$ be an orthonormal basis for $\mathcal{H}$. Then the set $B = \{(e_1,0), (e_2,0), \dots \} \cup \{(0,e_1), (0,e_2), \dots \}$ is an orthonormal basis for $\mathcal{H} \oplus \mathcal{H}$. Denote the complexification of $\mathcal{H} \oplus \mathcal{H}$ by $\widehat{\mathcal{H} \oplus \mathcal{H}}$, that is, $$\widehat{\mathcal{H} \oplus \mathcal{H}} = (\mathcal{H} \oplus \mathcal{H}) + i (\mathcal{H} \oplus \mathcal{H}).$$
    It it is easy to see that the set $B$ considered as a subset of $\widehat{\mathcal{H} \oplus \mathcal{H}}$ is an orthonormal basis for $\widehat{\mathcal{H} \oplus \mathcal{H}}$. For $n \in \mathbb{N}$, define the space $( \widehat{\mathcal{H} \oplus \mathcal{H}} )_{2n}$ as 
	$$( \widehat{\mathcal{H} \oplus \mathcal{H}} )_{2n} = \textrm{ span }\{(e_1, 0), (0, e_1), \dots, (e_n,0), (0,e_n) \}.$$ Let $P_{2n}$ denote the orthogonal projection of the Hilbert space $\widehat{\mathcal{H} \oplus \mathcal{H}}$ onto the subspace $( \widehat{\mathcal{H} \oplus \mathcal{H}} )_{2n}$. 
    Now
	$$
		i\hat{J} \leq I \Rightarrow \widehat{\sqrt{T}}i\hat{J}\widehat{\sqrt{T}} \leq \hat{T} \Rightarrow i\widehat{\sqrt{T}}\hat{J}\widehat{\sqrt{T}} \leq \hat{T} \Rightarrow P_{2n}(i\widehat{\sqrt{T}J\sqrt{T}})P_{2n} \leq P_{2n}\hat{T}P_{2n},$$
		that is, $(i\widehat{\sqrt{T}J\sqrt{T}})_{2n} \leq \hat{T}_{2n},$
	 where $\hat{T}_{2n}$ and $(i\widehat{\sqrt{T}J\sqrt{T}})_{2n}$ are the compression operators $P_{2n}\hat{T}P_{2n}$ and $P_{2n}(i\widehat{\sqrt{T}J\sqrt{T}})P_{2n}$, respectively, on $( \widehat{\mathcal{H} \oplus \mathcal{H}} )_{2n}$. These are operators on $2n$-dimensional spaces and hence can be identified with matrices of order $2n$. Now by using the Weyl's Monotonicity Theorem (Corollary III.2.3 of \cite{rbbook}), we have 
	\begin{equation} \label{weylineq}                \lambda_j^\downarrow((i\widehat{\sqrt{T}J\sqrt{T}})_{2n}) \leq \lambda_j^\downarrow(\hat{T}_{2n}), \quad j = 1,2, \dots, 2n,
	\end{equation} where $\lambda_j^\downarrow(A)$ denotes the eigenvalues of the matrix $A$ arranged in decreasing order.
	
	Now we make a more specific observation regarding the cardinality of the symplectic eigenvalues of $T$ lying strictly to the right (or left) of $\alpha$ with respect to the cardinality of the eigenvalues of $T$ lying strictly to the right (or left) of $\alpha$. Let $v_R,v_L$ respectively denote the number of eigenvalues of $T$ lying strictly to the right and left of $\alpha$; $s_R,s_L$ respectively denote the number of symplectic eigenvalues of $T$ lying strictly to the right and left of $\alpha$, where $v_R,v_L,s_R,s_L \in \{0,1,2,\dots\} \cup \{\infty\}$. 
	
	\begin{lemma}
		Let $\mathcal{H}$ be a real separable Hilbert space and $T$ be a positive invertible operator on $\mathcal{H} \oplus \mathcal{H}$ such that $T-\alpha I$ is compact for some $
		\alpha>0 $. Then
  \[s_R \leq v_R \quad \text{ and } \quad s_L \leq v_L,\]
  that is, the cardinality of the symplectic eigenvalues lying to the right (left) of $\alpha$ is at most the cardinality of usual eigenvalues lying to the right (left) of $\alpha$.
	\end{lemma}
	
	\begin{proof}
		If $v_R, v_L = \infty$, the statement is obvious. 
   Now suppose $v_R < \infty$. Enumerate the eigenvalues as follows
		$${\lambda_1^R}^\downarrow(T) \geq {\lambda_2^R}^\downarrow(T) \geq \dots \geq {\lambda_{v_R}^R}^\downarrow(T).$$
		Since they are finite in number, it follows from Theorem \ref{mnnresult} that for each $k = 1,2, \dots,$  as $n \rightarrow \infty$ the sequence of eigenvalues $\{ {\lambda_k^R}^\downarrow(T_n)\}$ converges to $\alpha$. Now applying limits on both sides of Equation \eqref{weylineq}, from Theorem \ref{mnnresult} we have for $j = 1,2, \dots, 	$
        \begin{equation} \label{eqinlemma}
            {\lambda_j^R}^\downarrow(i\widehat{\sqrt{T}J\sqrt{T}}) \leq {\lambda_j^R}^\downarrow(\hat{T}).
        \end{equation}

		But the symplectic eigenvalues of $T$ are the positive eigenvalues of the operator $i\widehat{\sqrt{T}J\sqrt{T}}$ (see Lemma \ref{lemmasymspecdefn}). Also, the eigenvalues of an operator and its complexification are the same. Hence, from Equation \eqref{eqinlemma}, we have 
		\begin{align*}
			{d_1^R}^\downarrow(T) &\leq {\lambda_1^{R}}^\downarrow(T) \\
			{d_2^R}^\downarrow(T) &\leq {\lambda_2^{R}}^\downarrow(T) \\
			&\vdots \\
			{d_{v_R}^R}^\downarrow(T) &\leq {\lambda_{v_R}^{R}}^\downarrow(T) \\
			{d_{v_R+1}^R}^\downarrow(T) &\leq \alpha \\
			{d_{v_R+2}^R}^\downarrow(T) &\leq \alpha \\
			&\vdots
		\end{align*} 
		Here we have considered the symplectic eigenvalues to the right of $\alpha$, so $\alpha \leq {d_j^R}^\downarrow(T)$ for all $j = 1,2, \dots.$ Hence, from the above set of inequalities, we have ${d_k^R}^\downarrow(T) = \alpha$ for all $k > v_R$, that is, the symplectic eigenvalues of $T$ lying strictly to the right of $\alpha$ is also finite and is at most $v_R$. Therefore $s_R < \infty$ when $v_R < \infty$ and $s_R \leq v_R$.
		
   Now suppose $v_L < \infty$. Enumerate the eigenvalues as follows
		$${\lambda_1^L}^\uparrow(T) \leq {\lambda_2^L}^\uparrow(T) \leq \cdots \leq {\lambda_{v_L}^L}^\uparrow(T).$$ Then the eigenvalues $(\lambda_j^L(T))^{-1} = \lambda_j^L(T^{-1})$ lie to the right of $\alpha^{-1}$ and we have $({\lambda_j^L}^\uparrow(T))^{-1} = {\lambda_j^R}^\downarrow(T^{-1})$. Note that 
  \begin{align*}          (i\widehat{\sqrt{T}J\sqrt{T}})^{-1} &= -i \widehat{(\sqrt{T}J\sqrt{T})^{-1}} = -i \widehat{((\sqrt{T})^{-1}J^{-1}(\sqrt{T})^{-1})} \\
  &= i \widehat{((\sqrt{T})^{-1}J(\sqrt{T})^{-1})}.
  \end{align*}
 Hence, $(d_j(T))^{-1} = d_j(T^{-1})$ and the symplectic eigenvalues of $T$ to the left of $\alpha$ are the symplectic eigenvalues of $T^{-1}$ to the right of $\alpha^{-1}$. Therefore we have $({d_j^L}^\uparrow(T))^{-1} = {d_j^R}^\downarrow(T^{-1})$. Now from the above part with $T^{-1}$ in the place of $T$, we have
		
		\begin{align*}
			({d_1^L}^\uparrow(T))^{-1} = {d_1^R}^\downarrow(T^{-1}) &\leq {\lambda_1^{R}}^\downarrow(T^{-1}) = ({\lambda_1^L}^\uparrow(T))^{-1}\\
			({d_2^L}^\uparrow(T))^{-1} = {d_2^R}^\downarrow(T^{-1}) &\leq {\lambda_2^{R}}^\downarrow(T^{-1}) = ({\lambda_2^L}^\uparrow(T))^{-1} \\
			&\vdots \\
			({d_{v_L}^L}^\uparrow(T))^{-1} = {d_{v_L}^R}^\downarrow(T^{-1}) &\leq {\lambda_{v_L}^R}^\downarrow(T^{-1}) = ({\lambda_{v_L}^L}^\uparrow(T))^{-1} \\
			({d_{{v_L}+1}^L}^\uparrow(T))^{-1} = {d_{{v_L}+1}^R}^\downarrow(T^{-1}) &\leq \alpha^{-1} \\
			({d_{{v_L}+2}^L}^\uparrow(T))^{-1} = {d_{{v_L}+2}^R}^\downarrow(T^{-1}) &\leq \alpha^{-1} \\
			&\vdots
		\end{align*} 
		
		But the symplectic eigenvalues ${d_j^R}^\downarrow(T^{-1})$ are bounded below by $\alpha^{-1}$, that is, $\alpha^{-1} \leq {d_j^R}^\downarrow(T^{-1})$ for all $j = 1,2, \dots$. Hence, from the set of inequalities above, we have $({d_k^L}^\uparrow(T))^{-1} = {d_k^R}^\downarrow(T^{-1}) = \alpha^{-1}$ for all $k > {v_L}$. That is, ${d_k^L}^\uparrow(T) = \alpha$ for all $k>{v_L}$, that is, the symplectic eigenvalues strictly to the left of $\alpha$ is also finite and is at most $v_L$. Therefore $s_L < \infty$ when $v_L < \infty$ and $s_L \leq v_L$. 
	\end{proof}
	
	Now we give some examples to illustrate  that strict inequality may hold in the previous lemma, that is, we show that the cardinality of the eigenvalues and the symplectic eigenvalues can be different.
 
	\begin{example} \label{exampleccfc} ($s_R < v_R = \infty, v_L = s_L = \infty$ ) \\
		Consider the operator $T$ defined on $\mathcal{H} \oplus \mathcal{H}$ (where $\mathcal{H}$ is an infinite-dimensional real separable Hilbert space) given by the matrix representation (with respect to some orthonormal basis) 
		$$T = 
		\begin{bmatrix}
			A & 0 \\
			0 & B
		\end{bmatrix},$$ where $A$ and $B$ are bounded real linear operators on $\mathcal{H}$ given by
		\begin{align*}
			A &= \textrm{ diag }\left\{3 + \frac{1}{n}: n \in \mathbb{N} \right\}, \\
			B &= \textrm{ diag }\left\{3 - \frac{1}{n}: n \in \mathbb{N} \right\}.
		\end{align*} Here $\alpha = 3$. Also for $j = 1,2, \dots,$ ${\lambda_j^R}^\downarrow(T) = 3 + \frac{1}{j}, {\lambda_j^L}^\uparrow(T) = 3 - \frac{1}{j}.$ That is, there are countably infinite eigenvalues to the right and left of $\alpha = 3$. The operator $i\widehat{\sqrt{T}J\sqrt{T}}$ is given by the matrix representation
		$$\begin{bmatrix}
			0 & i\widehat{A^\frac{1}{2}B^\frac{1}{2}} \\
			-i\widehat{B^\frac{1}{2}A^\frac{1}{2}} & 0
		\end{bmatrix} = \begin{bmatrix}
			0 & i\widehat{A^\frac{1}{2}B^\frac{1}{2}} \\
			-i\widehat{A^\frac{1}{2}B^\frac{1}{2}} & 0
		\end{bmatrix}$$ (since $A$ and $B$ commutes being diagonal operators).  The eigenvalues of the operator $i\widehat{\sqrt{T}J\sqrt{T}}$ are $$\left\{\pm \sqrt{\left(3 - \frac{1}{j}\right)\left(3 + \frac{1}{j}\right)}: j = 1,2, \dots \right\}.$$ Hence, the symplectic eigenvalues of $T$ are 
        $$\left\{ \sqrt{\left(3 - \frac{1}{j}\right)\left(3 + \frac{1}{j}\right)}: j = 1,2, \dots \right\}.$$ Since
        $\sqrt{\left(3 - \frac{1}{j}\right)\left(3 + \frac{1}{j}\right)} < 3$ for all $j = 1,2, \dots,$ all the symplectic eigenvalues of $T$ lie to the left of $3$. 
	\end{example}

	\begin{example} \label{exampleffff} ($s_R< v_R < \infty, v_L = s_L < \infty$) 
	\\ Consider the operator $T$ defined on $\mathcal{H} \oplus \mathcal{H}$ (where $\mathcal{H}$ is an infinite-dimensional real separable Hilbert space) given by the matrix representation (with respect to some orthonormal basis)
	$$
	T = \begin{bmatrix}
		A & 0 \\
		0 & B
	\end{bmatrix},
	$$ where $A$ and $B$ are bounded real linear operators represented by the infinite diagonal matrices given below
	$$
	A = 
	\begin{cases}
		\textrm{ diag } \{3 + \frac{1}{j}: j = 1,2, \dots, n \} \\
		\textrm{ diag } \{3: j> n \} 
	\end{cases}
	$$
	$$
	B = 
	\begin{cases}
		\textrm{ diag } \{3 - \frac{1}{j}: j = 1,2, \dots, m \} \\
		\textrm{ diag } \{3: j> m \} 
	\end{cases}
	$$ where $n$ and $m$ are finite positive integers such that $n > m$. Here $\alpha = 3$,
    $${\lambda_j^R}^\downarrow(T) = 3 + \frac{1}{j}, \quad j = 1,2, \dots, n$$
    and 
    $${\lambda_j^L}^\uparrow(T) = 3 - \frac{1}{j}, \quad j = 1,2, \dots, m.$$
    That is, there are finite eigenvalues to the right and left of $\alpha = 3.$ The operator $i\widehat{\sqrt{T}J\sqrt{T}}$ is given by the matrix representation
		$$\begin{bmatrix}
			0 & i\widehat{A^\frac{1}{2}B^\frac{1}{2}} \\
			-i\widehat{B^\frac{1}{2}A^\frac{1}{2}} & 0
		\end{bmatrix} = \begin{bmatrix}
			0 & i\widehat{A^\frac{1}{2}B^\frac{1}{2}} \\
			-i\widehat{A^\frac{1}{2}B^\frac{1}{2}} & 0
		\end{bmatrix}$$ (since $A$ and $B$ commutes being diagonal operators). The eigenvalues of the operator $i\widehat{\sqrt{T}J\sqrt{T}}$ are 
    \begin{align*}
        &\left\{\pm \sqrt{\left( 3 + \frac{1}{j}\right)\left( 3 - \frac{1}{j} \right)}, j = 1,2, \dots, m \right\} \\
        & \phantom{..........................................................} \cup \left\{ \pm \sqrt{3 \left(3 + \frac{1}{j} \right) }, j = m+1, \dots, n \right\} \cup \{3\}.
    \end{align*} 
    Hence, the symplectic eigenvalues of $T$ are 
    \begin{align*}
        &\left\{\sqrt{\left( 3 + \frac{1}{j}\right)\left( 3 - \frac{1}{j} \right)}, \quad j = 1,2, \dots, m \right\} \\
        &\phantom{...........................................................} \cup \left\{\sqrt{3 \left(3 + \frac{1}{j} \right) }, \quad j = m+1, \dots, n \right\} \cup \{3\}.
    \end{align*}   
    Note that 
    $$\sqrt{\left( 3 + \frac{1}{j}\right)\left( 3 - \frac{1}{j} \right)} < 3 \,\, \text{   and   } \,\,
    \sqrt{3 \left(3 + \frac{1}{j} \right)} > 3.$$ Hence, there are $m$ symplectic eigenvalues for $T$ strictly to the left of $\alpha = 3$ and $n - m$ symplectic eigenvalues for $T$ strictly to the right of $\alpha = 3$.
	\end{example}
	
	Now we state the main result of this article. Recall that $v_R$ and $v_L$ denote the number of eigenvalues of $T$ lying to the right and left of $\alpha$, respectively, and $s_R$ and $s_L$ denote the number of symplectic eigenvalues of $T$ lying to the right and left of $\alpha$, respectively.
	
	\begin{theorem} \label{decincgen}
		Let $\mathcal{H}$ be a real separable Hilbert space, $T$ be a positive invertible operator on $\mathcal{H} \oplus \mathcal{H}$ such that $T - \alpha I$ is compact for some $\alpha > 0$. Then
		\begin{enumerate}
			\item ${d_j^R}^\downarrow(T) \leq {\lambda_j^R}^\downarrow(T), \quad j = 1,2, \dots, s_R$; \label{decgen}
			\item ${\lambda_j^L}^\uparrow(T) \leq {d_j^L}^\uparrow(T), \quad j = 1,2, \dots, s_L.$ \label{incgen}
		\end{enumerate}
	\end{theorem}
	
	\begin{proof}
		Since $T - \alpha I$ is compact, $T$ has either finite eigenvalues or countably infinite eigenvalues.
		Now we prove the theorem by considering different possibilities as the following cases.
		
		\begin{enumerate}
			\item \underline{Case 1:} $v_R, v_L = \infty.$ 
			\\ Then the following cases can occur for the symplectic eigenvalues.
			\begin{enumerate}
				\item $s_R, s_L = \infty.$
				
				\item $s_R < \infty, s_L = \infty.$
				
				\item $s_R = \infty, s_L < \infty.$
				
				\item $s_R, s_L < \infty.$
			\end{enumerate}  Now we consider each of these sub-cases separately. 
			
			 \underline{For Case (a):} It is given that $s_R, s_L = \infty$, that is, $T$ has countably infinite symplectic eigenvalues strictly to the right and left of $\alpha$. By applying limits on both sides of Equation \eqref{weylineq}, from Theorem \ref{mnnresult} we have for $j = 1,2, \dots,$
			\begin{equation} \label{eqinequality}
				{\lambda_j^R}^\downarrow(i\widehat{\sqrt{T}J\sqrt{T}}) \leq {\lambda_j^R}^\downarrow(\hat{T})
			\end{equation}			
			From Lemma \ref{lemmasymspecdefn}, $\lambda_j(i \widehat{\sqrt{T}J\sqrt{T}})$ are the symplectic eigenvalues of $T$ lying symmetrically about $0$ on the real axis. So, ${\lambda_j^R}^\downarrow(i \widehat{\sqrt{T}J\sqrt{T}})$ are the symplectic eigenvalues of $T$. Also the eigenspectrum of an operator and its complexification are the same. Hence, from \eqref{eqinequality} we have $${d_j^R}^\downarrow(T) \leq {\lambda_j^R}^\downarrow(T), \quad j = 1,2, \dots.$$ 
			Now consider the eigenvalues and the symplectic eigenvalues of $T$ lying strictly to the left of $\alpha$. Since $T - \alpha I$ is compact, we have that the operator $$T^{-1} - \alpha^{-1}I = T^{-1} - \alpha^{-1}T^{-1}T 
			= T^{-1}(I - \alpha^{-1}T) 
			= -\alpha^{-1}T^{-1}(T - \alpha I)$$ is compact. Hence, it follows from Lemma \ref{lemmacomplexification} that $\widehat{T^{-1}}- \alpha^{-1} I$ is compact. Therefore $\widehat{T^{-1}}$ has a countable spectrum with $\sigma_{ess}(\widehat{T^{-1}}) = \{ \alpha^{-1} \}$. Then the eigenvalues $(\lambda_j^L(\hat{T}))^{-1} = \lambda_j^L(\widehat{T^{-1}})$  lies to the right of $\alpha^{-1}$ and we have $({\lambda_j^L}^\uparrow(\hat{T}))^{-1} = {\lambda_j^R}^\downarrow(\widehat{T^{-1}})$, where $j = 1,2, \dots$. Note that \begin{align*}          (i\widehat{\sqrt{T}J\sqrt{T}})^{-1} &= -i \widehat{(\sqrt{T}J\sqrt{T})^{-1}} = -i \widehat{((\sqrt{T})^{-1}J^{-1}(\sqrt{T})^{-1})} \\
  &= i \widehat{((\sqrt{T})^{-1}J(\sqrt{T})^{-1})}.
  \end{align*}
  Hence, $(d_j(T))^{-1} = d_j(T^{-1})$ and the symplectic eigenvalues of $T^{-1}$ lie to the right of $\alpha^{-1}$ with $({d_j^L}^\uparrow(T))^{-1} = {d_j^R}^\downarrow(T^{-1})$. Now from the above part with $T^{-1}$ in the place of $T$, we have
			$${d_j^R}^\downarrow(T^{-1}) \leq {\lambda_j^R}^\downarrow(T^{-1}), \qquad j = 1,2, \dots.$$ Therefore, 
			\begin{align*}
				({d_j^L}^\uparrow(T))^{-1} &\leq ({\lambda_{j}^L}^\uparrow(T))^{-1} \\
				\textrm{that is, } {\lambda_{j}^L}^\uparrow(T) &\leq {d_j^L}^\uparrow(T), \qquad j = 1,2, \dots.
			\end{align*}
			
			\underline{For Case(b):} Here it is assumed that $s_R < \infty, s_L = \infty$, that is, $T$ has finite symplectic eigenvalues strictly to the right of $\alpha$ and countably infinite symplectic eigenvalues strictly to the left of $\alpha$. Then the symplectic eigenvalues can be enumerated as $${d_1^R}^\downarrow(T) \geq {d_2^R}^\downarrow(T) \geq \cdots \geq {d_{s_R}^R}^\downarrow(T)$$ and is bounded below by $\alpha$. Now as in Case 1, by applying limits on both sides of Equation \eqref{weylineq}, we have for $j = 1,2, \dots,$
			$${\lambda_j^R}^\downarrow(i\widehat{\sqrt{T}J\sqrt{T}}) \leq {\lambda_j^R}^\downarrow(\hat{T}) = {\lambda_j^R}^\downarrow(T).$$ But from Theorem \ref{mnnresult}, the first $s_R$ terms of ${\lambda_j^R}^\downarrow(i \widehat{\sqrt{T}J\sqrt{T}})$ are ${d_j^R}^\downarrow(T)$ and the rest is $\alpha$. The inequality is still valid as the eigenvalues of $T$ is bounded below by $\alpha$. Hence, we have
			$$
			{d_j^R}^\downarrow(T) \leq {\lambda_j^R}^\downarrow(T) \quad \textrm{ for } j = 1,2, \dots, s_R.
			$$
			Now for the eigenvalues and the symplectic eigenvalues of $T$ which lie strictly to the left of $\alpha$, the relation is developed just as the corresponding formulations in Case 1(a).
			
			 \underline{For Case (c):} Here it is assumed that $s_R = \infty, s_L < \infty$, that is, $T$ has countably infinite symplectic eigenvalues strictly to the right of $\alpha$ and finite eigenvalues strictly to the left of $\alpha$. The proof is the same as in Case 1(b) with appropriate substitutions.
			
			\underline{For Case (d):} Here it is assumed that $s_R, s_L < \infty$, that is, $T$ has finite symplectic eigenvalues strictly to the right and left of $\alpha$. The relation follows by using the corresponding arguments from Cases 1(b) and 1(c). 
			
			\item \underline{Case 2:} $v_R < \infty, v_L = \infty$.
			\\ Then the following cases can occur for the symplectic eigenvalues.
			\begin{enumerate}
				
				\item $s_R < \infty, s_L = \infty$.  
				
				\item $s_R,s_L < \infty$.  
			\end{enumerate}
			Now we consider each of these sub-cases separately.
			
    \underline{For Case (a):} Here it is assumed that $s_R < \infty, s_L = \infty$, that is, $T$ has finite symplectic eigenvalues strictly to the right of $\alpha$ and countably infinite symplectic eigenvalues strictly to the left of $\alpha$. The relation between the eigenvalues and the symplectic eigenvalues of $T$ strictly to the left of $\alpha$ follows as in Case 1(a). Now for the eigenvalues and the  symplectic eigenvalues lying strictly to the right of $\alpha$, enumerate the eigenvalues to the right of $\alpha$ by
			$${\lambda_1^R}^\downarrow(T) \geq {\lambda_2^R}^\downarrow(T) \geq \cdots \geq {\lambda_{v_R}^R}^\downarrow(T).$$
			Now it follows from Theorem \ref{mnnresult} that for each $k = 1,2, \dots,$  as $n \rightarrow \infty$ the sequence of eigenvalues $\{ {\lambda_k^R}^\downarrow(T_n)\}$ and $\{ {\lambda_k^R}^\downarrow((i \widehat{\sqrt{T}J\sqrt{T}})_n) \}$ converges to $\alpha$. Now applying limits on both sides of Equation \eqref{weylineq}, from Theorem \ref{mnnresult} we have
			$${\lambda_j^R}^\downarrow(i\widehat{\sqrt{T}J\sqrt{T}}) \leq {\lambda_j^R}^\downarrow(\hat{T}) = {\lambda_j^R}^\downarrow(T), \qquad j = 1,2, \dots.$$ But since the eigenvalues and the symplectic eigenvalues strictly to the right of $\alpha$ are finite and the sequences $\{ {\lambda_k^R}^\downarrow(T_n)\}$ and $\{ {\lambda_k^R}^\downarrow((i \widehat{\sqrt{T}J\sqrt{T}})_n) \}$ converges to $\alpha$, we have ${\lambda_j^R}^\downarrow(T) = \alpha$ for all $j > v_R$ and 
			$$
			{\lambda_j^R}^\downarrow(i \widehat{\sqrt{T}J\sqrt{T}}) =
			\begin{cases}
				{d_j^R}^\downarrow(T), \quad j = 1,2, \dots, s_R; \\
				\alpha, \quad \,\, \qquad j > s_R,
			\end{cases}
			$$
			that is, 
			\begin{align*}
				{d_1^R}^\downarrow(T) &\leq {\lambda_1^{R}}^\downarrow(T) \\
				{d_2^R}^\downarrow(T) &\leq {\lambda_2^{R}}^\downarrow(T) \\
				&\vdots \\
				{d_{s_R}^R}^\downarrow(T) &\leq {\lambda_{s_R}^R}^\downarrow(T) \\
				\alpha &\leq {\lambda_{s_R + 1}^R}^\downarrow(T) \\
				&\vdots \\
				\alpha &\leq {\lambda_{v_R}^R}^\downarrow(T).
			\end{align*} 
			That is, $${d_j^R}^\downarrow(T) \leq {\lambda_j^R}^\downarrow(T), \quad j = 1,2, \dots, s_R.$$
			
   \underline{For Case(b):} Here it is given that $s_R,s_L < \infty$, that is, $T$ has finite eigenvalues strictly to the right and left of $\alpha$. The relation between the eigenvalues and the symplectic eigenvalues to the right of $\alpha$ follows just as in Case 2(a) and that of the eigenvalues and the symplectic eigenvalues to the left of $\alpha$ follows as in Case 1(c). 
			
			\item \underline{Case 3:} $v_R = \infty, v_L < \infty$. \\ Then the following cases can occur for the symplectic eigenvalues. 
			\begin{enumerate}
				\item $s_R = \infty, s_L < \infty$. 
				
				\item $s_R, s_L < \infty$.  
			\end{enumerate}
			
			\noindent Here the relations follow just as in Case 2 with appropriate substitutions.
			
			\item \underline{Case 4:} $v_R, v_L < \infty$. \\
			Here the only possible case is $s_R,s_L < \infty$, that is, $T$ has finite symplectic eigenvalues to the right and to the left of $\alpha$. Then the relation follows as in Case 2(a) for the eigenvalues and the symplectic eigenvalues to the right of $\alpha$ and as in Case 3(a) for the eigenvalues and the symplectic eigenvalues to the left of $\alpha$.
		\end{enumerate}
	\end{proof}

	\begin{remark} \label{finitedimcounterpart}
		Our result (Theorem \ref{decincgen}) when restricted to the finite-dimensional situation is equivalent to the  relation due to Bhatia and Jain \cite[Theorem 11(ii)]{rjt01}. This can be seen as follows. \newline The interlacing relation for matrices states that for any positive matrix $A$ of order $2n$, 
		\begin{equation} \label{interlacing}
			\lambda_j^\uparrow(A) \leq d_j^\uparrow(A) \leq \lambda_{n+j}^\uparrow(A), \quad j = 1,2, \dots, n
		\end{equation} where$\,\,d_1(A), d_2(A), \cdots, d_n(A)$ are the symplectic eigenvalues and $\lambda_1(A), \lambda_2(A), \dots, \lambda_{2n}(A)$ are the eigenvalues of $A$. Now the first set of inequalities in Equation \eqref{interlacing} gives 
		\begin{equation} \label{firstineq}
			\lambda_j^\uparrow(A) \leq d_j^\uparrow(A), \quad j = 1,2, \dots, n.
		\end{equation} Also the second set of inequalities in Equation \eqref{interlacing} gives $$d_j^\uparrow(A) \leq \lambda_{n+j}^\uparrow(A), \quad j = 1,2, \dots, n.$$ Now by arranging the eigenvalues and the symplectic eigenvalues of the equation above in decreasing order and rearranging the indices, we have 
		\begin{equation} \label{secondineq}
			d_j^\downarrow(A) \leq \lambda_j^\downarrow(A), \quad j = 1,2,\cdots, n
		\end{equation} 
        Choose $\alpha \in \mathbb{R}$ such that $\lambda_n^\uparrow(A) \leq \alpha \leq \lambda_{n+1}^\uparrow(A) (= \lambda_n^\downarrow(A))$. Rename the eigenvalues to the left of $\alpha$, $\lambda_j^\uparrow(A)$, as ${\lambda_j^L}^\uparrow(A)$ and to the right of $\alpha$, $\lambda_j^\downarrow(A)$, as ${\lambda_j^R}^\downarrow(A)$, where $j = 1,2, \dots , n$. Then Equations \eqref{firstineq} and \eqref{secondineq} becomes
        $${\lambda_j^L}^\uparrow(A) \leq d_j^\uparrow(A); \qquad d_j^\downarrow(A) \leq {\lambda_j^R}^\downarrow(A),$$ where $j = 1,2,\dots,n$.
	\end{remark}

	\section{Examples}\label{sec:special-cases}
In this section,  we illustrate our result on two important sub-classes of operators; the class of Gaussian Covariance Operators and the class of positive $\mathcal{AN}$ operators.

	\begin{definition} (Gaussian Covariance Operator \cite{bvtcrs01})  \label{defgco}
		Let $S$ be a real linear, bounded, symmetric and invertible operator on $\mathcal{H} \oplus \mathcal{H}$, where $\mathcal{H}$ is a real separable Hilbert space. Then $S$ is called a \emph{Gaussian Covariance Operator} (GCO) if the following three conditions are satisfied.
		\begin{enumerate}
			\item $\hat{S} - i\hat{J} \geq 0$, where $\hat{S}$, $\hat{J}$ are the complexification of the operators $S$ and $J$ respectively (See Definition \ref{defncomplexification}). 
			\item $S - I$ is Hilbert-Schmidt.
			\item $(\sqrt{S}J\sqrt{S})^T(\sqrt{S}J\sqrt{S}) - I$ is of trace class.
		\end{enumerate}
	\end{definition}
	\begin{remark} \label{cond3}
		The third condition in Definition \ref{defgco} can be replaced by the condition that $(JS)^2 + I$ is of trace class \cite[Corollary 3.3.1]{john2018infinite}.
	\end{remark}
	
	The Gaussian covariance operators are the covariance operators associated with quantum Gaussian states on a Hilbert space; see \cite{bvtcrs01} for a helpful characterization. The symplectic spectrum of the covariance operator	forms a complete invariant for Gaussian states as any two Gaussian states with the same symplectic spectrum conjugate with each other through a Gaussian symmetry \cite{bvtcrs01}.

	
	Lemma 3.3.2 of \cite{john2018infinite} states that the symplectic spectrum of any GCO lies to the right of $1$. Hence, the inequality in Theorem \ref{decincgen} corresponding to the symplectic eigenvalues lying to the left of $\alpha = 1$ is vacuously true in this case. Thus, we get the following reduction  by taking $T=S$ and $\alpha = 1$ in Theorem \ref{decincgen}.
	\begin{corollary} \label{decinc}
		For any Gaussian covariance operator $S$ on $\mathcal{H} \oplus \mathcal{H}$, its symplectic spectrum  lies to the right of $1$ and  $${d_j^R}^\downarrow(S) \leq {\lambda_j^R}^\downarrow(S), \quad j = 1, 2, \dots , s_R,$$
		where $s_R$ denotes the number of symplectic eigenvalues of $T$ lying strictly to the right $1$ and $s_R \in \{1,2, \dots \} \cup \{\infty \}.$
	\end{corollary}

	Absolutely Norm attaining operators ($\mathcal{AN}$ Operators) form an important class of infinite-dimensional operators. Here we show that invertible positive $\mathcal{AN}$ operators belong to the class of operators that we consider in this article. Hence, our main theorem applies to these operators as well. 

	\begin{definition} ($\mathcal{AN}$ Operators \cite{carvajal2012operators}) \label{defAN}
		Let $\mathcal{M}$ and $\mathcal{N}$ be complex Hilbert spaces. An operator $P \in \mathcal{B}(\mathcal{M}, \mathcal{N})$ is said to be an \emph{$\mathcal{AN}$ operator} or to satisfy the property $\mathcal{AN}$, if for every non-trivial closed subspace $\mathcal{E}$ of $\mathcal{M}$, there exists an element $x$ in $\mathcal{E}$ with unit norm such that $\|P\mid_{\mathcal{E}}\| = \|P\mid_{\mathcal{E}}(x)\|$.
	\end{definition}
	
	We are interested in positive $\mathcal{AN}$ operators ($(\mathcal{AN})_+$ operators). They form a proper cone in the real Banach space of Hermitian operators. A spectral characterization for $(\mathcal{AN})_+$ operators was formulated in \cite{pandey2017spectral} which we state below.
	
	\begin{theorem} \cite[Theorem 5.1]{pandey2017spectral}  \label{anoperatorsform}
		Let $\mathcal{H}$ be a complex Hilbert space of arbitrary dimension and let $P$ be a positive operator on $\mathcal{H}$. Then $P$ is an $\mathcal{AN}$ operator if and only if $P$ is of the form $P = \beta I + K + F$, where $\beta \geq 0$, $K$ is a positive compact operator and $F$ is a self-adjoint finite rank operator.
	\end{theorem}

    Williamson's normal form is defined for positive invertible operators on real separable Hilbert spaces. $(\mathcal{AN})_+$ operators are defined on complex Hilbert spaces. So to apply the normal form, we proceed as follows. Consider the complex Hilbert space $\mathcal{H}$. If $\mathcal{H}_\mathbb{R}$ is defined as the closure of the real span of $\mathcal{H}$, then $\mathcal{H} \simeq \mathcal{H}_\mathbb{R} \oplus \mathcal{H}_\mathbb{R}$ as real Hilbert spaces. We consider the positive invertible operator $P$ on the complex Hilbert space $\mathcal{H}$ as the real linear operator on the real Hilbert space $\mathcal{H}_{\mathbb{R}} \oplus \mathcal{H}_{\mathbb{R}}$. The same identification holds for the operators $K, F$ and $I$ in Theorem \ref{anoperatorsform}.

    From Theorem \ref{anoperatorsform} it is clear that if $\beta = 0$, then $P = K + F$ is compact and hence not invertible. Therefore, an $(\mathcal{AN})_+$ operator need not be invertible in general. Since the Williamson's normal form demands invertibility, we consider invertible $(\mathcal{AN})_+$ operators. That is, $\beta$ cannot be zero. Hence, $\beta > 0$. Then the operator $P$ takes the form $$P = \beta I + K + F, \quad \beta > 0.$$
    Since $K$ is a positive compact operator and $F$ is a finite rank self-adjoint operator, the sum $K^\prime = K + F$ is compact. That is, $P = K^\prime + \beta I, \beta > 0$. Now the result follows from Theorem \ref{decincgen} with $P$ in the place of $T$ and $\alpha = \beta.$

	\section{Concluding Remarks and Future Problems}
	
	\noindent Williamson's normal form and symplectic spectrum have found their importance in various fields of physics, quantum information theory being the latest interest. This article establishes a relationship between the symplectic eigenvalues and the eigenvalues for operators in a particular class, which contains the class of GCOs and $(\mathcal{AN})_+$ operators. This relationship is an infinite-dimensional analogue of Theorem 11(ii) in \cite{rjt01} and we notice interesting differences to the finite-dimensional result because of the difficulty in ordering the eigenvalues when we have infinitely many of them. Now we list down some future problems to consider.
	\begin{enumerate}
		
		\item Here we proved the interlacing relations for positive invertible operators $T$ on $\mathcal{H} \oplus \mathcal{H}$ such that $T - \alpha I$ is compact for some $\alpha > 0$. The next question is whether the relation holds for a much more general class of operators with a countable spectrum. Also, can we establish some bounds for operators with uncountable spectrum \cite{vbkant}. 
		
		\item Notice that in Theorem \ref{decincgen} we used the truncation method. This method is applicable to an arbitrary bounded self-adjoint operator. Since the operator $i\widehat{\sqrt{T}J\sqrt{T}}$ is a bounded self-adjoint operator, one can approximate the eigenvalues lying outside the bounds of the essential spectrum of $i\widehat{\sqrt{T}J\sqrt{T}}$. This indicates the scope for further applications of the truncation method in the symplectic spectral theory. 
		
		\item In this article, we developed an infinite-dimensional version of Theorem 11(ii) of \cite{rjt01}. The article \cite{rjt01} contains several other results for symplectic eigenvalues of matrices.  Finding analogues of these results to the infinite-dimension is an interesting project to consider.
		
		
	\end{enumerate} Furthermore, there is a scope for studying infinite-dimensional problems arising from the symplectic spectrum.

	

	




	\backmatter
	
	\bmhead{Acknowledgements}
	Tiju Cherian John thanks the Fulbright Scholar Program and United States-India Educational Foundation for providing funding and other support to conduct this research through a Fulbright-Nehru Postdoctoral Fellowship (Grant number: 2594/FNPDR/2020), he also acknowledges the United States Army Research Office MURI award on Quantum Network Science, awarded under grant number W911NF2110325 for partially funding this research. V. B. Kiran Kumar wishes to thank the SERB SURE Scheme (Project No. SUR/2022/003340) for the financial support. Anmary Tonny is supported by the INSPIRE PhD Fellowship of the Department of Science and Technology, Govt. of India. The authors would like to thank Akhila N. S. for her valuable suggestions. The authors are grateful to the Referee for his feedback, which helped us enhance the manuscript.
	
	\bmhead{Declarations}
	\begin{itemize}
		\item Conflict of interest/Competing interests: There is no competing interest.
		\item Availability of data and materials: There is no associated data.
	\end{itemize}

	\begin{appendices}
	\section{Inclusion Relations between the different Classes of Operators considered in this Article} \label{appendixinclusion}
	
	We have already remarked in Section \ref{discussionmainresults} that the inclusions of the operators that we consider are strict. Here we exhibit the same. 
 
 Let $\mathcal{H}$ be a real separable Hilbert space and $A$ be a positive invertible operator on $\mathcal{H} \oplus \mathcal{H}$ that has a matrix representation (with respect to some orthonormal basis) of the form
	$$A = \begin{bmatrix}
		B & 0 \\
		0 & B
	\end{bmatrix},$$ where the operator $B$ on $\mathcal{H}$ is given by the real matrix
	$$
	B = \textrm{ diag } \left\{ 3 - \frac{1}{n^2} : n \in \mathbb{N} \right\}.
	$$ Then by taking $\alpha = 3$, the operator $A$ comes under the bigger class but fails to be a GCO (this is because condition (2) of Definition \ref{defgco} fails) and an $(\mathcal{AN})_+$ operator (by Theorem \ref{anoperatorsform}). 
	
  Also if $A$ is taken as the operator with matrix representation (with respect to some orthonormal basis) 
	$$A = \begin{bmatrix}
		B & 0 \\
		0 & B
	\end{bmatrix},$$ where the operator $B$ on $\mathcal{H}$ is given by the real matrix
	$$
    B = \textrm{ diag } \left\{ 3 + \frac{1}{n^2} : n \in \mathbb{N} \right\},
	$$ then by taking $\alpha = 3$, $A$ becomes an $(\mathcal{AN})_+$ operator while it fails to be a GCO (condition (2) of Definition \ref{defgco} fails).
	Now if $A$ is taken as $$
	A = \begin{bmatrix}
		B & 0 \\
		0 & C
	\end{bmatrix},$$ where $B$ and $C$ are operators represented by the infinite real diagonal matrices \begin{align*}
		B &= \textrm{ diag } \left\{ 1 - \frac{1}{4^n} : n \in \mathbb{N} \right\}; \\
		C &= \textrm{ diag } \left\{ 1 + \frac{1}{2^n} : n \in \mathbb{N} \right\},
	\end{align*} then $\alpha = 1$ makes $A$ a GCO while it fails to be an $(\mathcal{AN})_+$ operator.
	
  Now if the operator $A$ is taken as the operator with matrix representation 
	$$A = \begin{bmatrix}
		B & 0 \\
    	0 & B
	\end{bmatrix},$$ where the operator $B$ is given by the real matrix
	$$
	B = \textrm{ diag } \left\{ 1 + \frac{1}{n^2} : n \in \mathbb{N} \right\},
	$$ then by taking $\alpha = 1$, $A$ is both a GCO and an $(\mathcal{AN})_+$ operator. 
 
 Let $G^\prime$ denote the subclass of the class of GCOs such that $G - I$ is positive for all $G \in G^\prime$. It is worth noticing that even though GCOs and $(\mathcal{AN})_+$ have no proper inclusions, $G^\prime$ is contained in the closure of $(\mathcal{AN})_+$ operators (this follows from Theorem 4.5 of \cite{doi:10.1080/03081087.2022.2126426} and the fact that the essential spectrum of a GCO is $\{1\}$). In fact, if we define a subclass $T^\prime$ of the class of operators we considered here such that $T - \alpha I$ is positive for all $T \in T^\prime$ and $\alpha> 0$, then $T^\prime$ is contained in the closure of $(\mathcal{AN})_+$ operator. 
	\end{appendices}




\begin{thebibliography}{1000}
		
		\bibitem{AdRaSa14}Adesso, G., Ragy, S., Lee, A. Continuous variable quantum information: Gaussian states and beyond. {\em Open Syst. Inf. Dyn.}. \textbf{21}, 1440001, 47 (2014), 
		https://doi.org/10.1142/S1230161214400010
		
		\bibitem{bv01}Bhat, B. V. R., John, T. C. Real normal operators and Williamson's normal form. {\em Acta Sci. Math. (Szeged)}. \textbf{85}, 507-518 (2019), https://doi.org/10.14232/actasm-018-570-5
		
		\bibitem{bvtcrs01}Bhat, B. V. R., John, T. C., Srinivasan, R. Infinite mode quantum Gaussian states. {\em Rev. Math. Phys.}. \textbf{31}, 1950030, 33 (2019), https://doi.org/10.1142/S0129055X19500302
		
		\bibitem{bhatiapositivebook} Bhatia, R. Positive Matrices. (Hindustan Book Agency TRIM 44. (2007)).
		
		\bibitem{rbbook} Bhatia, R. Matrix Analysis. (Springer Science \& Business Media, Vol 169. (2013))
		
		\bibitem{rjt01}Bhatia, R., Jain, T. On symplectic eigenvalues of positive definite matrices. {\em J. Math. Phys.}. \textbf{56}, 112201, 16 (2015), https://doi.org/10.1063/1.4935852
		
		
		
		\bibitem{bhatia_jain_2021}Bhatia, R., Jain, T. Variational principles for symplectic eigenvalues. {\em Canadian Mathematical Bulletin}. \textbf{64}, 553-559 (2021)
		
		\bibitem{bhatia2021szegHo}Bhatia, R., Jain, T., Sengupta, R. A Szegő type theorem and distribution of symplectic eigenvalues. {\em J. Spectr. Theory}. \textbf{11}, 1369-1389 (2021), https://doi.org/10.4171/jst/377
		
		\bibitem{botcher}Böttcher, A., Chithra, A., Namboodiri, M. N. N.  Approximation of approximation numbers by truncation. {\em Integral Equations Operator Theory}. \textbf{39}, 387-395 (2001), https://doi.org/10.1007/BF01203320
		
		
		\bibitem{carvajal2012operators}Carvajal, X., Neves, W. Operators that achieve the norm. {\em Integral Equations Operator Theory}. \textbf{72}, 179-195 (2012), https://doi.org/10.1007/s00020-011-1923-y
		
		\bibitem{gohberg} Gohberg, I., Goldberg, S., Kaashoek, M.A. Classes of Linear Operators (Vol. I, Birkhauser Verlag, Basel (1990))
		
		\bibitem{mag01}Gosson, M. Symplectic methods in harmonic analysis and in mathematical physics. (Birkhäuser/Springer Basel AG, Basel,2011), https://doi.org/10.1007/978-3-7643-9992-4
		
		\bibitem{logmaj}Hiai, F., Lim, Y. Log-majorizations for the (symplectic) eigenvalues of the Cartan barycenter. {\em Linear Algebra Appl.}. \textbf{553} pp. 129-144 (2018), https://doi.org/10.1016/j.laa.2018.04.029
		
		\bibitem{wolf}Idel, M., Soto Gaona, S., Wolf, M. Perturbation bounds for Williamson's symplectic normal form. {\em Linear Algebra Appl.}. \textbf{525} pp. 45-58 (2017), https://doi.org/10.1016/j.laa.2017.03.013
		
		
		\bibitem{jain2022derivatives}Jain, T., Mishra, H. K. Derivatives of symplectic eigenvalues and a Lidskii type theorem. {\em Canadian Journal Of Mathematics}. \textbf{74}, 457-485 (2022)
		
		\bibitem{john2018infinite}John, T. C. Infinite Mode Quantum Gaussian States. (Ph.D Thesis, Indian Statistical Institute-Kolkata, 2018)

        \bibitem{vbkant} Kiran Kumar, V. B., Tonny, Anmary. On approximating the symplectic spectrum of infinite-dimensional operators. \emph{J. Math. Phy.}. \textbf{65}, 042201 (2024), https://doi.org/10.1063/5.0169600 
		
		\bibitem{pandey2017spectral}Pandey, S., Paulsen, V. A spectral characterization of $\mathcal{AN}$ operators. {\em J. Aust. Math. Soc.}. \textbf{102}, 369-391 (2017), https://doi.org/10.1017/S1446788716000239
		
		\bibitem{Par13}Parthasarathy, K. R. The symmetry group of Gaussian states in $L^2(\mathbb{R}^n)$. {\em Prokhorov And Contemporary Probability Theory}. \textbf{33} pp. 349-369 (2013), http://dx.doi.org/10.1007/978-3-642-33549-521
		
		
		\bibitem{doi:10.1080/03081087.2022.2126426}Ramesh, G., Sequeira, S. S. On the closure of absolutely norm attaining operators. {\em Linear And Multilinear Algebra}. pp. 1-21 (2022), 
		
		
		\bibitem{tracemin}Son, N., Stykel, T. Symplectic eigenvalues of positive-semidefinite matrices and the trace minimization theorem. {\em Electron. J. Linear Algebra}. \textbf{38} pp. 607-616 (2022), https://doi.org/10.13001/ela.2022.7351
		
		
		
		\bibitem{weedbrooketal2012}Weedbrook, Christian and Pirandola, Stefano and Garc\'{\i}a Patr\'on, Ra\'ul and Cerf, Nicolas J. and Ralph, Timothy C. and Shapiro, Jeffrey H., Lloyd, Seth, Gaussian quantum information. {\em Rev. Mod. Phys.}. \textbf{84}, 621-669 (2012,5), https://link.aps.org/doi/10.1103/RevModPhys.84.621		
		\bibitem{jw01}Williamson, J. On the Algebraic Problem Concerning the Normal Forms of Linear Dynamical Systems. {\em Amer. J. Math.}. \textbf{58}, 141-163 (1936), https://doi.org/10.2307/2371062
		
		
		
	\end{thebibliography}
\end{document}